\documentclass[12pt]{article}

\usepackage{amsmath,amssymb,epsfig,fullpage}

\newtheorem{theorem}{Theorem}[section]

\newtheorem{corollary}[theorem]{Corollary}

\numberwithin{equation}{section}

\newcommand{\N}{\mathbb{N}}

\newcommand{\C}{\mathbb{C}}

\newcommand{\ov}{\overline}
\newcommand{\dis}{\displaystyle}
\newcommand{\supp}{\textup{supp}}

\newcommand{\Notequiv}{/\kern-.6em\hbox{$\equiv$} }
\newcommand{\db}{\rule[.05in]{.09in}{.10in}} 

\begin{document}

\title{Reverse Triangle Inequalities for Potentials}%
\author{I. E. Pritsker\thanks{Research partially supported
by the National Security Agency (grant H98230-06-1-0055), and by the
Alexander von Humboldt Foundation.}\ \ and E. B. Saff
\thanks{Research supported, in part, by the National Science Foundation (DMS-0603828).}\\ \\
{\it Dedicated to George G. Lorentz, whose works have been a great
inspiration}}%

\date{}%

\maketitle


\noindent{\bf Abstract.} We study the reverse triangle inequalities
for suprema of logarithmic potentials on compact sets of the plane.
This research is motivated by the inequalities for products of
supremum norms of polynomials. We find sharp additive constants in
the inequalities for potentials, and give applications of our
results to the generalized polynomials.

We also obtain sharp inequalities for products of norms of the
weighted polynomials $w^nP_n,\ \deg(P_n)\le n,$ and for sums of
suprema of potentials with external fields. An important part of our
work in the weighted case is a Riesz decomposition for the weighted
farthest-point distance function.


\section{Products of polynomials and sums of potentials} \label{sec1}

Let $E$ be a compact set in the complex plane ${\C}$. Given the
bounded above functions $f_j, j=1,\ldots,m,$ on $E$, we have by a
standard inequality that
\[
\sup_E \sum_{j=1}^m f_j \le \sum_{j=1}^m \sup_E f_j.
\]
It is not possible to reverse this inequality for arbitrary
functions, even if one introduces additive or multiplicative
``correction" constants. However, we are able to prove the reverse
inequalities for {\em logarithmic potentials}, with sharp additive
constants. For a positive Borel measure $\mu$ with compact support
in the plane, define its (subharmonic) potential \cite[p. 53]{Ra} by
\[
p(z):=\int \log|z-t|\,d\mu(t).
\]
Let $\nu_j, j=1,\ldots,m,$ be positive compactly supported Borel
measures with potentials $p_j$. We normalize the problem by assuming
that the total mass of $\nu:=\sum_{j=1}^m \nu_j$ is equal to 1, and
consider the inequality
\begin{align} \label{1.1}
\sum_{j=1}^m \sup_E p_j \le C_E(m) + \sup_E \sum_{j=1}^m p_j.
\end{align}
Clearly, if \eqref{1.1} holds true, then $C_E(m)\ge 0.$ One may also ask whether \eqref{1.1} holds with a constant $C_E$ independent of $m.$
The motivation for such inequalities comes directly from
inequalities for the norms of products of polynomials. Indeed, if
$P(z)=\prod_{j=1}^n (z-a_j)$ is a monic polynomial, then $\log
|P(z)| = n \int \log|z-t|\,d\tau(t)$. Here, $\tau=\frac{1}{n}
\sum_{j=1}^n \delta_{a_j}$ is the normalized counting measure in the
zeros of $P$, with $\delta_{a_j}$ being the unit point mass at
$a_j.$ Let $\|P\|_E:=\sup_E |P|$ be the uniform (sup) norm on $E$.
Thus \eqref{1.1} takes the following form for polynomials $P_j,
j=1,\ldots,m,$
\[
\prod_{j=1}^m \| P_j \|_E \leq e^{nC_E(m)} \left\| \prod_{j=1}^m
 P_j \right\|_E,
\]
where $n$ is the degree of the product polynomial $\prod_{j=1}^m
 P_j$. We outline a brief history of such inequalities below.

Kneser \cite{Kn} proved the first sharp inequality for norms of
products of polynomials on $[-1,1]$ (see also Aumann \cite{Au} for a
weaker result)
\begin{equation}  \label{1.2}
\| P_1 \|_{[-1,1]} \|P_2 \|_{[-1,1]} \leq K_{\ell ,n} \| P_1 P_2
\|_{[-1,1]}, \quad \deg P_1 = \ell,\ \deg P_2 =n-\ell,
\end{equation}
where
\begin{equation}  \label{1.3}
K_{\ell ,n} := 2^{n -1} \prod_{k =1}^{\ell} \left( 1 + \cos \frac{2k
-1}{2n} \pi \right) \prod_{k =1}^{n - \ell} \left( 1 + \cos \frac{2k
-1}{2n} \pi \right).
\end{equation}
Observe that equality holds in (\ref{1.2}) for the Chebyshev
polynomial $t(x) = \cos n \arccos x = P_1 (x) P_2 (x)$, with a
proper choice of the factors $P_1 (x)$ and $P_2 (x)$. Borwein
\cite{Bor} generalized this to the multifactor inequality
\begin{equation} \label{1.4}
\prod_{j =1}^m \| P_j \|_{[-1,1]} \leq 2^{n -1} \prod_{k =1}^{[
\frac{n}{2} ]} \left( 1 + \cos \frac{2k -1}{2n} \pi \right)^2
\left\| \prod_{j=1}^m  P_j \right\|_{[-1,1]},
\end{equation}
where $n$ is the degree of $\prod_{j=1}^m  P_j$. We remark that
\begin{equation} \label{1.5}
2^{n -1} \prod_{k =1}^{[ \frac{n}{2} ]} \left( 1 + \cos \frac{2k
-1}{2n} \pi \right)^2 \sim (3.20991 \ldots )^n \mbox{ as } n
\rightarrow \infty.
\end{equation}

Another inequality of this type for $E = D$, where $D := \{ w: |w|
\le 1 \}$ is the closed unit disk, was proved by Gelfond \cite[p.
135]{Ge} in connection with the theory of transcendental numbers:
\begin{equation} \label{1.6}
\prod_{j =1}^m \| P_j \|_{D} \leq e^n \left\| \prod_{j=1}^m
 P_j \right\|_{D}.
\end{equation}
Mahler \cite{Ma1} later replaced $e$ by $2$:
\begin{equation} \label{1.7}
\prod_{j =1}^m \| P_j \|_{D} \leq 2^n \left\| \prod_{j=1}^m
 P_j \right\|_{D} .
\end{equation}
It is easy to see that the base $2$ cannot be decreased, if $m =n$
and $n \rightarrow \infty$.  However, Kro\'{o} and Pritsker
\cite{KP} showed that, for any $m \leq n,$
\begin{equation} \label{1.10}
\prod_{j=1}^m \| P_j \|_{D} \leq 2^{n -1} \left\| \prod_{j=1}^m
 P_j \right\|_{D},
\end{equation}
where equality holds in (\ref{1.10}) for {\it each} $n \in {\N}$,
with $m=n$ and $\prod_{j=1}^m P_j = z^n -1$. On the other hand, Boyd
\cite{Boy1, Boy2} proved that, given the number of factors $m$ in
(\ref{1.7}), one has
\begin{equation} \label{1.8}
\prod_{j=1}^m \| P_j \|_{D} \leq (C_m)^n \left\| \prod_{j=1}^m
 P_j \right\|_{D},
\end{equation}
where
\begin{equation} \label{1.9}
C_m := \exp \left( \frac{m}{\pi} \int_0^{\pi/m} \log \left(2 \cos
\frac{t}{2}\right) dt \right)
\end{equation}
is asymptotically best possible for {\em each fixed} $m$, as $n
\rightarrow \infty$.

For a compact set $E\subset\C,$ a natural general problem is to
find the {\it smallest} constant $M_E > 0$ such that
\begin{equation} \label{1.11}
\prod_{j=1}^m \| P_j \|_E \leq (M_E)^n \left\| \prod_{j=1}^m
 P_j \right\|_E
\end{equation}
holds for arbitrary polynomials $\{ P_j (z) \}_{j=1}^m$ with complex
coefficients, where $n = \deg(\prod_{j=1}^m P_j)$. The solution of
this problem is based on the logarithmic potential theory (cf.
\cite{Ra} and \cite{Ts}). Let ${\rm cap}(E)$ be the {\it logarithmic
capacity} of a compact set $E \subset {\C}$. For $E$ with ${\rm
cap}(E)>0$, denote the {\it equilibrium measure} of $E$  by $\mu_E$.
We remark that $\mu_E$ is a positive unit Borel measure supported on
the outer boundary of $E$ (see \cite[p. 55]{Ts}). Define
\begin{equation} \label{1.12}
d_E(z) := \max_{t \in E} |z -t|, \qquad z \in {\C},
\end{equation}
which is clearly a positive and continuous function in ${\C}$. It is
easy to see that the logarithm of this distance function is
subharmonic in $\C.$ Moreover, it has the following integral
representation
\[
\log d_E(z) = \int \log |z -t| d \sigma_E(t), \quad z \in {\C} ,
\]
where $\sigma_E$ is a positive unit Borel measure in ${\C}$ with
unbounded  support, see Lemma 5.1 of \cite{Pr1} and \cite{LP01}.
Further study of the representing measure $\sigma_E$ is contained in
the work of Gardiner and Netuka \cite{GN}. This integral
representation is the key fact used by the first author to prove the
following result.

\begin{theorem} \label{thm1.1} \textup{\cite{Pr1}}
Let $E \subset {\C}$ be a compact set, ${\rm cap} (E) >0$.  Then the
best constant $M_E$ in {\rm (\ref{1.11})} is given by
\begin{equation} \label{1.13}
M_E = \frac{\exp\left(\dis\int \log d_E(z) d \mu_E (z)\right)}{{\rm
cap} (E)} .
\end{equation}
\end{theorem}
It is not difficult to see that $M_E$ is invariant under the
similarity transformations of the set $E$ in the plane.

For the closed unit disk $D$, we have that ${\rm cap}(D) =1$ and
that $d\mu_{D} = d\theta/(2 \pi),$ where $d\theta$ is the arclength
on $\partial D$ \cite[p. 84]{Ts}.  Thus Theorem \ref{thm1.1} yields
\begin{equation*}
M_{D} = \exp\left(\frac{1}{2 \pi} \int_0^{2 \pi} \log d_{D}(e^{i
\theta})\ d \theta\right) = \exp\left(\frac{1}{2 \pi} \int_0^{2 \pi}
\log 2\ d \theta\right) =2,
\end{equation*}
so that we immediately obtain Mahler's inequality (\ref{1.7}).

If $E = [-1,1]$ then ${\rm cap}([-1,1]) = 1/2$ and $d\mu_{[-1,1]} =
dx/(\pi \sqrt{1 -x^2}) , \ x \in [-1,1],$ which is the Chebyshev
distribution (see \cite[p. 84]{Ts}). Using Theorem \ref{thm1.1}, we
obtain
\begin{eqnarray*}
M_{[-1,1]} & = & 2\exp\left(\frac{1}{\pi} \int_{-1}^1 \frac{\log
d_{[-1,1]}(x)}{\sqrt{1 -x^2}} dx\right) = 2\exp\left(\frac{2}{\pi}
\int_0^1 \frac{\log (1 +x)}{\sqrt{1 -x^2}} dx \right) \nonumber \\
& = & 2\exp\left(\frac{2}{\pi} \int_0^{\pi/2} \log (1 + \sin t) dt
\right) \approx 3.2099123,
\end{eqnarray*}
which gives the asymptotic version of Borwein's inequality
(\ref{1.4})-(\ref{1.5}).

Considering the above analysis of Theorem \ref{thm1.1}, it is
natural to conjecture that the sharp universal bounds for $M_E$ are
given by
\begin{equation} \label{1.14}
2=M_{D} \le M_E \le M_{[-1,1]} \approx 3.2099123,
\end{equation}
for any bounded non-degenerate continuum $E$, see \cite{Pr3}. This
problem was treated in the recent papers of the first author and
Ruscheweyh \cite{PR1} and \cite{PR2}, where the lower bound $M_E \ge
M_D =2$ is proved for all compact sets $E$, and the upper bound is
proved for certain special classes of continua. A general approach
to this type of extremal problem was proposed by Baernstein,
Laugesen and Pritsker \cite{BLP}. We show in the next section that
all results about $M_E$ are directly applicable to the constants
$C_E$ and $C_E(m)$ in the inequality for potentials \eqref{1.1}.

The assumption that $E$ is of positive capacity is vital for our
results. For example, when $E$ consists of a finite number of points
$\{z_j\}_{j=1}^N,\ N\ge 2$, then no inequality of the type
\eqref{1.11} is possible with {\em any} constant. Indeed, if $m=n\ge
N$ then we consider $P_j(z)=z-z_j,\ j=1,\ldots,N,$ and $P_j(z)\equiv
1, j>N,$ which gives $\|P_j\|_E > 0,\ j=1,\ldots,m,$ but $\|
\prod_{j=1}^m P_j \|_E = 0.$

For infinite countable sets $E$ we have cap$(E)=0$, and the constants in the inequalities for norms of products of polynomials may grow arbitrarily fast.

\begin{theorem} \label{thm1.2}
Let $\{A_n\}_{n=1}^{\infty}$ be any increasing sequence satisfying $A_n\ge 1$. There exists an infinite countable set $E$ such that
\begin{eqnarray} \label{1.15}
\sup_{P_j} \frac{\prod_{j=1}^m \|P_j\|_E}{\|\prod_{j=1}^m P_j\|_E} \ge A_n,\quad n=\deg\left(\prod_{j=1}^m P_j\right)\in\N.
\end{eqnarray}
\end{theorem}

Thus one should expect faster-than-exponential growth of constants, if the assumption cap$(E)>0$ is lifted.

\section{Main results} \label{sec2}

Our first inequality stated in Theorem \ref{thm2.1} includes the constant $C_E$ that is independent of the number of potentials $m.$ In fact, Theorem \ref{thm2.1} may be deduced from our Theorem \ref{thm2.4}, which takes $m$ into account, and gives a sharp version of \eqref{1.1}.
\begin{theorem} \label{thm2.1}
Let $E \subset {\C}$ be a compact set, ${\rm cap}(E)>0$. Suppose
that $\nu_j, j=1,\ldots,m,$ are positive compactly supported Borel
measures with potentials $p_j$, such that the total mass of
$\sum_{j=1}^m \nu_j$ is equal to 1. We have
\begin{align} \label{2.1}
\sum_{j=1}^m \sup_E p_j \le C_E + \sup_E \sum_{j=1}^m p_j,
\end{align}
where
\begin{equation} \label{2.2}
C_E := \int \log d_E(z) d \mu_E (z) - \log{\rm cap}(E)
\end{equation}
cannot be replaced by a smaller constant.
\end{theorem}
Since $C_E$ is independent of $m$, it is possible to extend
\eqref{2.1} to infinite sums of potentials. One should ensure the
absolute convergence of the series $\sum_{j=1}^{\infty} p_j$ on $E$
for this purpose.

We note that $C_E$ is invariant under the similarity transforms of
the plane, i.e. under the maps $\phi(z)=az+b$ or $\phi(z)=a\bar{z}+b,$ where $a,b\in\C.$ It is obvious from \eqref{1.13} that $C_E = \log M_E$.
Hence the results of \cite{BLP,PR1,PR2} apply here, and we obtain
the following.

\begin{corollary} \label{cor2.2}
Let $E \subset {\C}$ be an arbitrary compact set, ${\rm cap} (E)
>0$. Then $C_E\ge \log 2,$ where equality holds if and only if
$\partial U \subset E \subset U$, where $U$ is a closed disk.
\end{corollary}

\begin{corollary} \label{cor2.3}
Let $E \subset {\C}$ be a connected compact set, but not a single
point. Suppose that $z,w\in E$ satisfy $\textup{diam}\,E=|z-w|$ and
the line segment $[z,w]$ joining $z$ to $w$ lies in $E$. If $E$ is
contained in the disk with diameter $[z,w]$, then
\[
C_E \le \frac{2}{\pi} \int_0^2 \frac{\log(2+x)}{\sqrt{4-x^2}}\,dx =
C_{[-2,2]} \approx \log{3.2099123}.
\]
Furthermore, this inequality holds for any centrally symmetric
continuum $E$ that contains its center of symmetry.
\end{corollary}

We conjecture in line with \eqref{1.14} (see \cite{Pr3,PR1}) that
$C_E \le C_{[-2,2]}$ for all non-degenerate continua $E$. One may readily see that $C_E \le \log (\textup{diam}\,E/\textup{cap}(E)) \le \log 4$ for any non-degenerate continuum. Improved estimates may also be found in \cite{BLP,PR1,PR2}.

We now explore the dependence of $C_E(m)$ in \eqref{1.1} on $m$. The
key results for a polynomial analog are due to Boyd \cite{Boy1,
Boy2} for the unit disk, see \eqref{1.8}-\eqref{1.9}. The polynomial
case for general sets was touched upon in \cite{Pr1}, and developed
further in \cite{PR2}.

A closed set $S\subset E$ is called {\em dominant} if
\begin{equation} \label{2.3}
d_E(z)=\max_{t\in S} |z-t| \quad\mbox{ for all } z \in
\textup{supp}(\mu_E).
\end{equation}
When $E$ has at least one finite dominant set, we define a {\em minimal dominant} set $\mathfrak D_E$ as a dominant set with the smallest number of points, i.e. card$(\mathfrak D_E)$ is minimal. Of course, $E$ might not have finite dominant sets at all, in which case we can take any
dominant set as the minimal dominant set with card$(\mathfrak
D_E)=\infty$, e.g., $\mathfrak D_E=\partial E.$

\begin{theorem} \label{thm2.4}
Let $E \subset {\C}$ be compact, ${\rm cap} (E) >0$. Suppose that
$\nu_j, j=1,\ldots,m,$ are positive compactly supported Borel
measures with potentials $p_j$, such that the total mass of
$\sum_{j=1}^m \nu_j$ is equal to 1. Then
\begin{align} \label{2.4}
\sum_{j=1}^m \sup_E p_j \le C_E(m) + \sup_E \sum_{j=1}^m p_j,
\end{align}
where
\begin{align} \label{2.5}
C_E(m):=\max_{c_k\in\partial E} \int \log \max_{1 \leq k \leq m} |z
-c_k|\, d \mu_E (z) - \log {\rm cap}(E)
\end{align}
cannot be replaced by a smaller constant for each fixed $m\ge 2.$
Furthermore, if $m < \textup{card}(\mathfrak D_E)$ then $C_E(m) <
C_E,$ while $C_E(m) = C_E$ for $m \ge \textup{card}(\mathfrak D_E)$.
When $\mathfrak D_E$ is infinite, $C_E(m) < C_E$ holds for all
$m\in\N,\ m\ge 2.$
\end{theorem}

Since $|z - c_k| \le d_E(z),\ c_k\in\partial E,$ it is clear from
\eqref{2.2} and \eqref{2.5} that $C_E(m) \le C_E$ for all $E$ and
all $m\in\N.$ Thus Theorem \ref{thm2.1} is an immediate consequence of
Theorem \ref{thm2.4}. If the sets $\{c_k\}_{k=1}^m$ are dense in $\partial
E$ as $m\to\infty$, then $\dis\lim_{m\to\infty} \max_{1 \leq k \leq
m} |z -c_k| = d_E(z),\ z\in\C.$ Hence $\dis\lim_{m\to\infty} C_E(m)
= C_E.$ However, the following result shows that we always have
strict inequality for smooth sets.

\begin{corollary} \label{cor2.5}
If $E\subset\C$ is a compact set bounded by finitely many
$C^1$-smooth curves, then $C_E(m) < C_E$ for all $m\in\N,\ m\ge 2.$
\end{corollary}

On the other hand, we have $C_E(m) = C_E$ for $m \ge s$ for every
polygon with $s$ vertices. Furthermore, not all vertices may belong
to the minimal dominating set. For example, if $E$ is an obtuse
triangle, then $\mathfrak D_E$ consists of only two vertices that
are the endpoints of the longest side. Hence $C_E(m) = C_E$ for $m
\ge 2$ as in the segment case. Any circular arc of angular measure
at most $\pi$ has its endpoints as the minimal dominating set, which
gives $C_E(m) = C_E$ for $m \ge 2$ here too. However, if the angular
measure of this arc is greater than $\pi$, then one immediately
obtains that $\mathfrak D_E$ is infinite, and $C_E(m) < C_E$ for all
$m\ge 2.$

Finding the exact values of $C_E(m)$ for general sets is very
complicated. It is analogous to finding solutions of discrete energy
problems. Following Boyd \cite{Boy1,Boy2}, we give the values of
$C_D(m)$, where $D$ is a disk, see \eqref{1.8}-\eqref{1.9}.

\begin{corollary} \label{cor2.6}
If $E$ is a closed disk $D$, then
\[
C_D(m)=\frac{m}{\pi} \int_0^{\pi/m} \log \left(2 \cos
\frac{t}{2}\right) dt, \quad m\ge2.
\]
\end{corollary}
It is easy to see that $C_D(m)<C_D=\log 2, \ m\ge2.$

We conclude this section with an application of our results for
potentials to generalized polynomials of the form
$P_j(z)=\prod_{k=1}^{k_j} |z-z_{k,j}|^{r_k},$ where $k_j\in\N$ and
$z_{k,j}\in\C,\ r_k>0,\ k=1,\ldots,k_j.$ Let $n_j:=\sum_{k=1}^{k_j}
r_k$ be the degree of the generalized polynomial $P_j.$

\begin{corollary} \label{cor2.7}
Let $E \subset {\C}$ be a compact set, ${\rm cap} (E) >0$. If $P_j,\ j=1,\ldots,m,$ are the generalized polynomials of the corresponding degrees $n_j,$ then
\begin{align*}
\prod_{j=1}^m \| P_j \|_E \le e^{nC_E(m)} \left\| \prod_{j=1}^m
 P_j \right\|_E \le e^{nC_E} \left\| \prod_{j=1}^m
 P_j \right\|_E,
\end{align*}
where $n=\sum_{j=1}^m n_j$, and where $C_E$ and $C_E(m)$ are defined by \eqref{2.2} and \eqref{2.5} respectively.
\end{corollary}
We remind the reader that $C_E=\log M_E$ , so that the above corollary extends Theorem \ref{thm1.1}.

\section{Weighted polynomials and potentials} \label{sec3}

In this section, we assume that $E\subset\C$ is any closed set,
which is not necessarily bounded. Let $w:E\to[0,\infty)$ be an {\em
admissible  weight function} \cite[p. 26]{ST} in the sense of
potential theory with external fields. This means that
\begin{itemize}
\item $w$ is upper semicontinuous on $E$
\item $\textup{cap}\left(\left\{z\in E : w(z)>0\right\}\right) > 0$
\item If $E$ is unbounded then $\dis\lim_{|z|\to\infty, z\in E} |z|w(z) = 0$
\end{itemize}
It is implicit that cap$(E)>0$ in this case. We study certain
analogs of our main results for weighted polynomials of the form
$w^k(z)P(z),$ $\deg(P)\le k,$ as well as for potentials with
external fields. In order to state such analogs, we need the notions
of the weighted equilibrium measure $\mu_w$ and the modified Robin's
constant $F_w$. Recall that $\mu_w$ is a positive unit Borel measure
supported on a compact set $S_w\subset E,$ that is characterized by the inequalities
\[
\int \log|z-t|\,d\mu_w(t) + \log w(z) + F_w \ge 0,\quad z\in
S_w=\supp(\mu_w),
\]
and
\[
\int \log|z-t|\,d\mu_w(t) + \log w(z) + F_w \le 0,\quad \mbox{for
q.e. }z\in E,
\]
where q.e. (quasi everywhere) means that the above inequality holds
with a possible exceptional set of zero capacity (cf. Theorem 1.3 of
\cite[p. 27]{ST}). We refer to \cite{ST} for a detailed survey of
potential theory with external fields. The {\em weighted farthest-point
distance function}
\begin{align} \label{3.0}
d_E^w(z):=\sup_{t\in E} w(t) |z-t|,\quad z\in\C,
\end{align}
plays an important role in our results, resembling the role of its
predecessor $d_E(z)$ defined in \eqref{1.12}.

\begin{theorem} \label{thm3.1}
Let $E \subset {\C}$ be a closed set, and let $w$ be an admissible
weight on $E$. If $P_j,\ j=1,\ldots,m,$ are polynomials of the
corresponding degrees $n_j,$ then
\begin{align} \label{3.1}
\prod_{j=1}^m \| w^{n_j} P_j \|_E \le e^{nC_E^w(m)} \left\| w^{n} \prod_{j=1}^m P_j \right\|_E \le e^{nC_E^w} \left\| w^{n} \prod_{j=1}^m
 P_j \right\|_E,
\end{align}
where $n=\sum_{j=1}^m n_j$. The constant
\begin{align} \label{3.2}
C_E^w(m):=\sup_{c_k\in E} \int \log \max_{1 \leq k \leq m} w(c_k)|z
-c_k|\, d \mu_w(z) + F_w
\end{align}
cannot be replaced by a smaller value for each fixed $m\ge 2.$
Also,
\begin{align} \label{3.3}
C_E^w:=\int \log d_E^w(z)\, d \mu_w(z) + F_w
\end{align}
cannot be replaced by a smaller value independent of $m.$
\end{theorem}
If $w$ is continuous and $E$ has positive capacity at each of its points,
then any weighted polynomial of the form $w^kP,$ $\deg(P)\le
k,$ attains its norm on $S_w$, which is often a {\em proper} subset
of $E$ (cf. \cite{MS} and Section III.2 of \cite{ST}). More
generally, the norm is always attained on $S_w\cup \ov R_w\subset
E,$ where $R_w:=\{z\in E: \int \log|z-t|\,d\mu_w(t) + \log w(z) + F_w
>0\}$, see Theorem 2.7 of \cite[p. 158]{ST}. Thus all sup norms in
Theorem \ref{thm3.1} may be replaced by the norms on $S_w\cup \ov
R_w$. As a consequence for the weighted distance function $d_E^w$,
we observe that for any $z\in\C$ there exists $\zeta_z\in S_w\cup
\ov R_w$ such that
\[
d_E^w(z) = \|w(\cdot)(z-\cdot)\|_E = \|w(\cdot)(z-\cdot)\|_{S_w\cup
\ov R_w} = w(\zeta_z)|z-\zeta_z|.
\]

We give a couple of examples of the constant $C_E^w$ for specific
sets and weights below. It is clear that any example of this kind
heavily depends on the knowledge of the weighted equilibrium measure
$\mu_w$ and the modified Robin's constant $F_w$. In addition, one
should be
able to compute the weighted distance function $d_E^w$.

\medskip
\noindent{\bf Examples\\
1. Incomplete polynomials of G. G. Lorentz:} Let $E=[0,1]$ and $w(x)=x$. It is known that
\[
d\mu_w(x) = \frac{2}{\pi x} \sqrt {\frac{x-1/4}{1-x}},\quad x\in
S_w=[1/4,1],
\]
see \cite[p. 243]{ST}. We also have that $F_w = 8\log2-3\log3$ by
\cite[p. 206]{ST}. Furthermore, it follows from a direct calculation
that
\[
d_E^w(x) = \left\{\begin{array}{ll} 1-x, \quad &1/4\le x\le
2(\sqrt{2}-1), \\
x^2/4,\quad &2(\sqrt{2}-1)\le x \le 1.\end{array} \right.
\]
The approximate numerical value obtained from \eqref{3.3} is $C_E^w
\approx 1.037550517,$ so that \eqref{3.1} gives
\[
\prod_{j=1}^m \| x^{n_j} P_j(x) \|_{[0,1]} \le (2.8222954)^{n}
\left\| x^{n} \prod_{j=1}^m  P_j(x) \right\|_{[0,1]},
\]
where $\deg P_j \leq n_j$ and $n=n_1+\ldots+n_m.$ (The polynomials
$x^{n_j}P_j(x)$ are special examples of incomplete polynomials, a
subject that was introduced by G. G. Lorentz in \cite{Lo}.)
Note that the above inequality is a significant improvement of the
Borwein-Kneser inequality \eqref{1.4}-\eqref{1.5} applied to the
polynomials $x^{n_j} P_j(x)$ on $[0,1].$ Indeed, since the degree of
$x^{n_j} P_j(x)$ equals $2n_j,$ we obtain from
\eqref{1.4}-\eqref{1.5} (or from \eqref{1.13}) that
\[
\prod_{j=1}^m \| x^{n_j} P_j(x) \|_{[0,1]} \le (10.303537)^{n}
\left\| x^{n} \prod_{j=1}^m  P_j(x) \right\|_{[0,1]},
\]
where the constant comes from $(M_{[0,1]})^2 \approx (3.2099123)^2 <
10.303537.$

\medskip
\noindent {\bf 2.} Let $E=\C$ and $w(z)=e^{-|z|}$. In this case, we have
\cite[p. 245]{ST} that
\[
d\mu_w(re^{i\theta}) = \frac{1}{2\pi}\,dr\,d\theta,\quad r\in[0,1],\
\theta\in[0,2\pi),
\]
$d_E^w(z)=e^{|z|-1}$ for $z\in S_w=\{z:|z|\le 1\}$, and $F_w=1.$
Here we explicitly find that $C_E^w = 1/2$ and consequently, from
\eqref{3.1},
\[
\prod_{j=1}^m \| e^{-n_j|z|} P_j(z) \|_{\C} \le e^{n/2} \left\|
e^{-n|z|} \prod_{j=1}^m P_j(z) \right\|_{\C}.
\]

\bigskip
With the notation of Theorems \ref{thm2.1}-\ref{thm2.4}, we let
$\alpha_j:=\nu_j(\C)$ be the total mass of the measure $\nu_j.$ For
the potentials with external fields $p_j(z)+\alpha_j \log w(z),$ we
have the following estimates.

\begin{theorem} \label{thm3.2}
Let $E \subset {\C}$ be a closed set, and let $w$ be an admissible weight on $E$. Suppose that $\nu_j, j=1,\ldots,m,$ are positive compactly supported Borel measures with potentials $p_j$, such that $\nu_j(\C)=\alpha_j$ and $\sum_{j=1}^m \alpha_j = 1.$ Then
\begin{align} \label{3.4}
\sum_{j=1}^m \sup_E (\alpha_j \log w + p_j) &\le C_E^w(m) + \sup_E \left( \log w + \sum_{j=1}^m p_j \right) \nonumber \\ &\le C_E^w + \sup_E \left( \log w + \sum_{j=1}^m p_j \right).
\end{align}
The constants $C_E^w$ and $C_E^w(m)$ are defined by \eqref{3.2} and
\eqref{3.3} respectively. They are sharp here in the same sense as
in Theorem \ref{thm3.1}.
\end{theorem}
Using a well known connection between polynomials and potentials of
discrete measures, we observe that Theorem \ref{thm3.1} is a direct
consequence of Theorem \ref{thm3.2}. For each polynomial
$P_j(z)=\prod_{k=1}^{n_j} (z-z_{k,j}),\ j=1,\ldots,m,$ we associate
the zero counting measure $\nu_j:=\frac{1}{n}\sum_{k=1}^{n_j}
\delta_{z_{k,j}}.$ Since
\[
\frac{1}{n}\log|P_j(z)| = \int\log|z-t|\,d\nu_j(t)=p_j(z) \quad
\mbox{and} \quad \frac{1}{n}\log\|w^{n_j}P_j\|_E = \sup_E
\left(\frac{n_j}{n}\log w + p_j\right),
\]
it is clear that \eqref{3.4} gives the log of \eqref{3.1} in this
notation. Another immediate observation is that Theorem \ref{thm3.2}
implies \eqref{2.1} and \eqref{2.4}, if we set $w\equiv 1$ on $E$.

The key ingredient in our proofs of Theorems \ref{thm3.1} and
\ref{thm3.2} is the following Riesz representation for $\log
d_E^w(z),$ which may be of independent interest.

\begin{theorem} \label{thm3.3}
Let $E \subset {\C}$ be a closed set. Suppose that $w:E\to[0,\infty)$ is upper semicontinuous on $E$, and that $w\not\equiv 0$ on $E$. If $E$ is unbounded then we also assume that $\dis\lim_{|z|\to\infty, z\in E} |z|w(z) = 0$. The function
\begin{align} \label{3.5}
\log d_E^w(z) := \sup_{t\in E} \left(\log w(t) + \log|z-t|\right) =
\log \|w(\cdot)(z-\cdot)\|_E,\quad z\in\C,
\end{align}
is subharmonic in $\C$, and
\begin{align} \label{3.6}
\log d_E^w(z)=\int \log |z-t|\, d\sigma_E^w(t) + \sup_E \log w,\quad z\in\C,
\end{align}
where $\sigma_E^w$ is a positive unit Borel measure.
\end{theorem}
Note that we relaxed conditions on the weight $w$ in Theorem
\ref{thm3.3} by not requiring the set $\left\{z\in E :
w(z)>0\right\}$ be of positive capacity. Such weights are called
quasi-admissible in \cite{ST}. Since the proofs of Theorems
\ref{thm3.1} and \ref{3.2} only require \eqref{3.6} for a finite set
$E=\{c_k\}_{k=1}^m$, we give a short and transparent proof of this
special case. The complete general proof of Theorem \ref{thm3.3}
will appear in our forthcoming work, together with a comprehensive
study of the weighted distance function.

We remark that $d_E^w$ is Lipschitz continuous in $\C$, which
readily follows from triangle inequality. Indeed, we have that
$|z_1-t| \le |z_1-z_2|+|z_2-t|$ for all $z_1,z_2\in\C$ and all $t\in
E.$ Hence
\[
d_E^w(z_1) = \sup_{t\in E} w(t)|z_1-t| \le |z_1-z_2|\sup_{t\in E}
w(t) + \sup_{t\in E} w(t)|z_2-t| = |z_1-z_2|\sup_{E} w + d_E^w(z_2)
\]
and
\[
|d_E^w(z_2) - d_E^w(z_1)| \le |z_2-z_1| \sup_{E} w,
\]
after interchanging $z_1$ and $z_2.$ If the set $\left\{z\in E :
w(z)>0\right\}$ is not a single point, then $d_E^w$ is strictly
positive in $\C,$ and $\log d_E^w$ is also Lipschitz continuous in
$\C.$ In particular, this always holds for admissible weights.

\section{Proofs} \label{secP}

\noindent{\it Proof of Theorem \ref{thm1.2}.} Without loss of generality we assume that $A_n\to\infty$ as $n\to\infty.$ Consider the sequence $x_1=1$ and $x_n=1/(2A_n),\ n\ge 2,$ and let $E:= \{x_n\}_{n=1}^{\infty} \cup \{0\}.$ Thus $E$ is a compact subset of $[0,1].$ Define $P_j(x):=x-x_j,\ j\in\N.$ Note that $\|P_1\|_E=1$, $\|P_j\|_E = 1-x_j \ge \frac{1}{2},\ j\ge 2,$ and
\[
\left\| \prod_{j=1}^n P_j \right\|_E = \prod_{j=1}^n x_j \le \frac{1}{2^{n-1} A_n},\quad n\in\N.
\]
Hence
\[
\frac{\prod_{j=1}^n \|P_j\|_E}{\left\|\prod_{j=1}^n P_j\right\|_E} = \frac{\prod_{j=2}^n (1-x_j)}{\prod_{j=1}^n x_j} \ge \frac{2^{-n+1}}{2^{-n+1} A_n^{-1}} = A_n.
\]

\db

\noindent{\it Proofs of Theorems \ref{thm2.1} and \ref{thm2.4}.}
Since any subharmonic potential $p_k$ is upper semicontinuous, it
attains a supremum on the compact set $E.$ Furthermore, we can
assume that the supremum is attained on $\partial E$, by the maximum
principle for subharmonic functions \cite[p. 29]{Ra}. Thus for any $k =1,
\ldots ,m$, there exists $c_k \in \partial E$ such that
\[
\sup_E p_k = p_k (c_k ).
\]
Define the function
\[
d_m (z) := \max_{1 \leq k \leq m} |z-c_k|, \qquad z \in {\C}.
\]
Lemma 2 of Boyd \cite{Boy2} states that for any set $\{ c_k \}_{k =1}^m\subset\C$ there exists a probability measure $\sigma_m$ such that
\[
\log d_m (z) = \int \log |z -t| d \sigma_m (t), \qquad z \in {\C}.
\]
Let $\nu:= \sum_{k =1}^m \nu_k,$ so that $\nu$ is a unit measure with the potential $p(z)=\int |z-t|\,d\nu(t)=\sum_{k =1}^m p_k(z).$ We
use the integral representation of $\log d_m$ and Fubini's theorem
in the following estimate:
\begin{eqnarray} \label{5.1}
\sum_{k =1}^m \sup_E p_k &=& \sum_{k=1}^m p_k (c_k) = \sum_{k =1}^m
\int \log |c_k-z|\,d\nu_k(z) \leq \int \log d_m (z)\,d\nu(z)
\nonumber \\
&=& \int \int \log |z-t|\,d\sigma_m(t)\,d\nu(z) = \int p(t)\,
d\sigma_m (t).
\end{eqnarray}
It is known \cite{Pr1} that the support of $\sigma_m$ is unbounded,
so that we need to estimate the growth of $p$ in the plane by its
supremum on the set $E$. This is analogous to the Bernstein-Walsh
lemma for polynomials \cite[p. 156]{Ra}. Let
$g(t):=\dis\int\log|t-z|\,d\mu_E(z) - \log{\rm cap}(E),\ t\in\C,$
and note that $g(t) \ge 0, \ t\in E,$  by Frostman's theorem
\cite[p. 59]{Ra}. On the other hand, we trivially have that $p(t) - \sup_E
p \le 0,\ t\in E,$ which gives
\[
g(t) \ge p(t) - \sup_E p,\ t\in E.
\]
By the Principle of Domination  (see \cite[p. 104]{ST}), we deduce
that the last inequality holds everywhere:
\[
p(t) \le \sup_E p + g(t),\ t\in\C.
\]
This inequality applied in \eqref{5.1} yields
\begin{eqnarray*}
\sum_{k =1}^m \sup_E p_k &\leq& \int \left(\sup_E p + g(t)\right)
d\sigma_m(t) = \sup_E p + \int g(t)\,d\sigma_m(t)\\ &=& \sup_E p +
\int\left(\int \log |z -t|\, d\mu_E (z) - \log {\rm cap}(E)\right) d\sigma_m(t)\\
&=& \sup_E p + \int \int \log|z -t|\,d\mu_E(z)\,d\sigma_m(t) - \log {\rm cap}(E)\\
&=& \sup_E p + \int \int \log|z -t|\,d\sigma_m(t)\,d\mu_E(z) - \log {\rm cap}(E) \\
&=& \sup_E p + \int \log d_m (z)d \mu_E (z) - \log {\rm cap}(E),
\end{eqnarray*}
where we consecutively used $\sigma_m(\C)=1,$ the representation of
$g$ via the potential of $\mu_E$, Fubini's theorem, and the integral
representation for $\log d_m$. Hence \eqref{2.4} follows from the
above estimate by taking maximum over all possible $m$-tuples of
$c_k\in\partial E,\ k=1,\ldots,m.$ (Note that $\log d_m$ is a
continuous function in the variables $c_k$, so that $\dis\int \log
d_m(z)d \mu_E(z)$ is continuous too.) Furthermore, \eqref{2.1} is
immediate after observing that $d_m(z)\le d_E(z),\ z\in\C.$

Suppose that $\dis\int \log d_m(z)d \mu_E(z)$ attains its maximum on
$(\partial E)^m$ for some set $c_k^*\in\partial E,\ k =1, \ldots
,m.$ We now show that $C_E(m)$ cannot be replaced by a smaller
constant for a fixed $m\ge 2.$ Let
\[
d_m^*(z):=\max_{1 \leq k \leq m} |z - c_k^*|, \quad z \in {\C},
\]
and define the sets
\[
S_1:=\{z\in\supp\,\mu_E : |z-c_1^*|=d_m^*(z) \}
\]
and
\[
S_k:=\{z\in\supp\,\mu_E\setminus \cup_{j=1}^{k-1} S_j :
|z-c_k^*|=d_m^*(z) \},\quad k=2,\ldots,m.
\]
It is clear that
\[
\supp\,\mu_E=\bigcup_{k=1}^m S_k \quad\mbox{and}\quad S_k \bigcap
S_j = \emptyset,\ k\neq j.
\]
Hence the measures $\nu_k^*:=\mu_E\vert_{S_k}$ give the
decomposition
\[
\mu_E = \sum_{k=1}^m \nu_k^*.
\]
If $E$ is regular, then $\dis\int\log|z-t|\,d\mu_E(z) = \log {\rm
cap}(E), \ t\in E,$ by Frostman's theorem \cite[p. 59]{Ra}. Thus we obtain
that
\begin{align*}
\sum_{k=1}^m \sup_E p_k^* &\ge \sum_{k=1}^m p_k^*(c_k^*) =
\sum_{k=1}^m \int \log|c_k^*-z|\,d\nu_k^*(z) = \sum_{k=1}^m \int
\log d_m^*(z)\,d\nu_k^*(z) \\ &=  \int \log d_m^*(z)\,d\mu_E(z) -
\log {\rm cap}(E) + \sup_{t\in E} \dis\int\log|z-t|\,d\mu_E(z).
\\ &=C_E(m) + \sup_E \sum_{k=1}^m p_k^*.
\end{align*}
Hence equality holds in \eqref{2.4} in this case.

An alternative proof that $C_E(m)$ cannot be replaced by a smaller
constant for any set $E$ (that does not require $E$ to be regular)
may be given by using the $n$-th Fekete points ${\mathcal F}_n = \{
a_{l,n} \}_{l =1}^n$ of $E$ \cite[p. 152]{Ra}. Let $\{c_k^*\}_{k=1}^m$ be
the maximizers of $\dis\int \log d_m(z)d \mu_E(z)$ on $(\partial
E)^m$, as before. We define a subset ${\mathcal F}_{k,n} \subset \{
a_{l,n} \}_{l =1}^n,$ associated with each point $c_k^*,\
k=1,\ldots,m,$ so that $a_{l,n} \in {\mathcal F}_{k,n}$ if
\begin{equation} \label{5.2}
d_m^*(a_{l,n})=|a_{l,n} - c_k^*|, \quad 1 \leq l \leq n.
\end{equation}
In the case that \eqref{5.2} holds for more than one $c_k^*$, we
assign $a_{l,n}$ to only one set ${\mathcal F}_{k,n}$, to avoid an
overlap of these sets. It is clear that, for any $n \in \N$,
\[
\bigcup_{k=1}^m {\mathcal F}_{k,n} = {\mathcal F}_n \quad {\rm and}
\quad {\mathcal F}_{k_1,n} \bigcap {\mathcal F}_{k_2,n} = \emptyset,
\ k_1 \neq k_2.
\]
Define the measures
\[
\nu_{k,n}^*:=\frac{1}{n} \sum_{a_{l,n} \in {\mathcal F}_{k,n}}
\delta_{a_{l,n}},
\]
so that for their potentials
\[
p_{k,n}^*(z) = \frac{1}{n} \sum_{a_{l,n} \in {\mathcal F}_{k,n}}
\log|z -a_{l,n}|, \quad k=1,\ldots,m,
\]
we have
\[
\sup_E p_{k,n}^* \ge \frac{1}{n} \sum_{a_{l,n} \in {\mathcal
F}_{k,n}} \log|c_k^* - a_{l,n}| = \frac{1}{n} \sum_{a_{l,n} \in
{\mathcal F}_{k,n}} \log d_m^*(a_{l,n}), \quad k=1,\dots,m.
\]
It follows from the weak* convergence of $\nu_n^*:=\sum_{k=1}^m
\nu_{k,n}^* = \frac{1}{n} \sum_{l=1}^n \delta_{a_{l,n}}$ to $\mu_E$,
as $n\to\infty,$ that
\begin{align*}
\liminf_{n \rightarrow \infty} \sum_{k=1}^m \sup_E p_{k,n}^* &\ge
\lim_{n \rightarrow \infty} \frac{1}{n} \sum_{k =1}^n \log d_m^*(a_{k,n}) \\
&= \int \log d_m^*(z) d \mu_E (z).
\end{align*}
Also, we have for the potential $p_n^*$ of $\nu_n^*$ that \cite[p. 155]{Ra}
\[
\lim_{n \rightarrow \infty} \sup_E p_n^* = \lim_{n \rightarrow
\infty} \log \left\| \prod_{k =1}^n (z-a_{k,n}) \right\|_E^{1/n} =
\log\textup{cap}(E),
\]
which gives
\[
\liminf_{n \rightarrow \infty} \sum_{k=1}^m \sup_E p_{k,n}^* \ge
C_E(m) + \lim_{n \rightarrow \infty} \sup_E p_n^*.
\]
Hence \eqref{2.4} turns into asymptotic equality as $n\to\infty,$
with $m\ge 2$ being fixed.

A similar argument with Fekete points shows that \eqref{2.1} turns
into asymptotic equality when $m=n\to\infty.$

Since $d_m (z) \leq d_E(z)$ for any $z \in {\C}$, we immediately
obtain that $C_E(m)\le C_E.$ Suppose that $m < \textup{card}
(\mathfrak D_E).$ Then there is $z_0\in\textup{supp}(\mu_E)$ such that
$d_m^* (z_0) < d_E(z_0).$  As both functions are continuous, the
same strict inequality holds in a neighborhood of $z_0$, so that
$\dis\int \log d_m^*(z)\, d \mu_E (z) < \dis\int \log d_E(z)\, d
\mu_E (z)$ and $C_E(m)<C_E.$ When $\mathfrak D_E$ is infinite, this
argument gives that $C_E(m)<C_E,\ m\ge 2.$ Assume now that
$\mathfrak D_E$ is finite and that $m \ge \textup{card} (\mathfrak
D_E).$ Then $d_m^*(z)=d_E(z)$ for all $z\in\textup{supp}(\mu_E),$
because one of the possible choices of the points $\{ c_k \}_{k
=1}^m\subset
\partial E$ includes points of the set $\mathfrak D_E$. It is
immediate that $\dis\int \log d_m^*(z)\, d \mu_E (z) = \dis\int \log
d_E(z)\, d \mu_E (z)$ and $C_E(m)=C_E$ in this case.

\db

\noindent{\it Proof of Corollary \ref{cor2.2}.} Use $C_E=\log M_E$
and apply Theorem 2.5 of \cite{PR1}.

\db

\noindent{\it Proof of Corollary \ref{cor2.3}.} The first part when
$E$ is contained in the disk with diameter $[z,w]$ follows from
$C_E=\log M_E$ and Corollary 2.2 of \cite{PR1}. The second part for
centrally symmetric $E$ is a consequence of Corollary 6.3 from
\cite{BLP}.

\db

\noindent{\it Proof of Corollary \ref{cor2.5}.} We need to show that the minimal dominant set is infinite, hence the result follows from
Theorem \ref{thm2.4}. Suppose to the contrary that $\mathfrak
D_E=\{\zeta_l\}_{l=1}^s$ is finite. Let $J\subset\partial E$ be a
smooth closed Jordan curve. Define
\[
J_l:=\{z\in J: d_E(z)=|z-\zeta_l|\}. \quad l=1,\ldots,s.
\]
It is clear that $J=\cup_{l=1}^s J_l.$ Observe that the segment
$[z,\zeta_l],\ z\in J_l,$ is orthogonal to $\partial E$ at
$\zeta_l$. Hence each $J_l$ is contained in the normal line to
$\partial E$ at $\zeta_l,\ l=1,\ldots,s.$ We thus obtain that $J$ is
contained in a union of $s$ straight lines, so that $J$ cannot have
a continuously turning tangent, which contradicts the smoothness
assumption.

\db

\noindent{\it Proof of Corollary \ref{cor2.6}.} Apply Theorem 1 of
\cite{Boy2}, and use that $C_D(m)=\log C_m$, where $C_m$ is given in
\eqref{1.9}.

\db

\noindent{\it Proof of Corollary \ref{cor2.7}.} For
$P_j(z)=\prod_{k=1}^{k_j} |z-z_{k,j}|^{r_k},$ define the zero
counting measures
\[
\nu_j = \frac{1}{n} \sum_{k=1}^{k_j} r_k \delta_{z_{k,j}},\quad
j=1,\ldots,m,
\]
where $\delta_z$ is a unit point mass at $z.$ We obtain that
\[
p_j(z) = \int \log|z-t|\,d\nu_j(t) = \frac{1}{n} \log |P_j(z)|,\quad j=1,\ldots,m.
\]
Hence Corollary \ref{cor2.7} follows by combining \eqref{2.1} and \eqref{2.4}.

\db

\noindent{\it Proofs of Theorems \ref{thm3.1} and \ref{thm3.2}.} We
follow some ideas used to prove Theorem \ref{thm1.1} in \cite{Pr1}
and Theorems \ref{thm2.1}-\ref{thm2.4} of this paper, augmented with
certain necessary facts on weighted potentials and distance
function. Note that admissibility of the weight $w$ implies
$\dis\lim_{z\to\infty,z\in E} \log w(z) - \log|z| = -\infty$.
Combining this observation with upper semicontinuity of $\log w$ and
of the potentials $p_j$, we conclude that there exist points $c_j \in E$
satisfying
\[
\sup_E (\alpha_j \log w + p_j) = \alpha_j \log w(c_j) +
p_j(c_j),\quad j=1,\ldots,m.
\]
Consider the weighted distance function
\[
d_m^w (z) := \max_{1 \leq j \leq m} w(c_j)|z-c_j|, \qquad z \in
{\C},
\]
and write by Theorem \ref{thm3.3} (with $E=\{c_j\}_{j=1}^m$ there) that
\begin{align} \label{5.3}
\log d_m^w (z) = \int \log |z -t| d \sigma_m^w (t) + \max_{1 \leq j
\leq m} \log w(c_j), \qquad z \in {\C},
\end{align}
where $\sigma_m^w$ is a probability measure on ${\C}$. Consider the
unit measure $\nu:= \sum_{k =1}^m \nu_k$ and its potential
$p(z)=\dis\int \log|z-t|\,d\nu(t)=\sum_{k =1}^m p_k(z).$ Using
\eqref{5.3} and Fubini's theorem, we have
\begin{align} \label{5.6}
\sum_{j=1}^m \sup_E (\alpha_j \log w + p_j) &= \sum_{j=1}^m
\left(\alpha_j \log w(c_j) + p_j(c_j)\right) \nonumber \\ &=
\sum_{j=1}^m \left(\alpha_j \log w(c_j) + \int \log
|c_j-z|\,d\nu_j(z)\right) \nonumber \\ &\leq \int \log
d_m^w(z)\,d\nu(z) \nonumber \\
&= \int \int \log |z-t|\,d\sigma_m^w(t)\,d\nu(z) + \max_{1 \leq j
\leq m} \log w(c_j) \nonumber \\ &= \int p(t)\, d\sigma_m^w (t) +
\max_{1 \leq j \leq m} \log w(c_j).
\end{align}
We now need an estimate of $p$ in $\C$ via the sup of $\log w + p$
on $E$. Obviously, $\log w(t) + p(t) \le \sup_E (\log w + p)$ for
$t\in S_w,$ as $S_w\subset E.$ We also know from Theorem 1.3 of
\cite[p. 27]{ST} that
\[
\dis\int\log|t-z|\,d\mu_w(z) +F_w \ge -\log w(t),\quad t\in S_w.
\]
This gives
\[
p(t) \le \sup_E (\log w + p) - \log w(t) \le \sup_E (\log w + p) +
\dis\int\log|t-z|\,d\mu_w(z) +F_w,\quad t\in S_w.
\]
Hence we have the desired estimate
\[
p(t) \le \sup_E (\log w + p) + \dis\int\log|t-z|\,d\mu_w(z)
+F_w,\quad t\in\C,
\]
by the Principle of Domination \cite[p. 104]{ST}. We proceed with
inserting the above inequality into \eqref{5.6}, and estimate as
follows
\begin{align*}
\sum_{j=1}^m \sup_E (\alpha_j \log w + p_j) &\leq \int \left(\sup_E
(\log w + p) + \dis\int\log|t-z|\,d\mu_w(z) + F_w\right) d\sigma_m^w(t) \\
&+ \max_{1 \leq j \leq m} \log w(c_j)\\
&= \sup_E (\log w + p) + F_w + \max_{1 \leq j \leq m} \log w(c_j) \\
&+ \int \int \log|z -t|\,d\sigma_m^w(t)\,d\mu_w(z)\\
&= \sup_E (\log w + p) + F_w + \int \log d_m^w(z)\,d\mu_w(z),
\end{align*}
where we again used $\sigma_m^w(\C)=\mu_w(\C)=1,$ the representation
for $\log d_m^w$, and Fubini's theorem. Hence the first inequality
in \eqref{3.4} follows by taking sup over $m$-tuples of $c_j\in E,\
j=1,\ldots,m.$ The second inequality is immediate from $d_m^w(z)\le
d_E^w(z),\ z\in\C.$

It was explained after the statement of Theorem \ref{thm3.2} that
Theorem \ref{thm3.1} is its special case. In particular, we have
that \eqref{3.4} for the zero counting measures $\nu_j$ of
polynomials $P_j$ implies \eqref{3.1}. Thus \eqref{3.1} is also
proved. On the other hand, if we show that the constants $C_E^w(m)$
and $C_E^w$ are sharp in Theorem \ref{thm3.1}, then they are
obviously sharp in Theorem \ref{thm3.2} too. Hence we select this
path and prove sharpness for the weighted polynomial case, i.e., for
discrete measures in weighted Fekete points.

Since $\log d_m^w(z)$ is an upper semicontinuous function of $c_j\in
E,\ j=1, \ldots,m,$ we have that $\dis\int \log d_m^w(z)\, d
\mu_E(z)$ is also upper semicontinuous in those variables, and hence
attains its maximum on $E^m$ for some set $c_j^*\in E,\
j=1,\ldots,m.$ We now show that $C_E(m)$ cannot be replaced by a
smaller constant for each fixed $m$, by adapting the proof of
Theorem 4.1 from \cite{Pr1}. Let
\[
d_m^*(z):=\max_{1 \leq j \leq m} w(c_j^*)|z - c_j^*|, \qquad z \in {\C}.
\]
Consider the weighted $n$-th Fekete points ${\mathcal F}_n=\{ a_{l,n} \}_{l=1}^n$ for the weight $w$ on $E$, and the corresponding polynomials (cf. Section III.1 of \cite{ST})
\[
F_n(z) = \prod_{l=1}^n (z-a_{l,n} ), \quad n \in {\N}.
\]
We define a subset ${\mathcal F}_{j,n} \subset \{ a_{l,n} \}_{l
=1}^n$ associated with each point $c_j^*,\ j=1,\ldots,m,$ so that
$a_{l,n} \in {\mathcal F}_{j,n}$ if
\begin{equation} \label{5.5}
d_m^*(a_{l,n})=w(c_j^*) |a_{l,n} - c_j^*|, \quad 1 \leq l \leq n.
\end{equation}
If \eqref{5.5} holds for more than one $c_j^*$, then we include
$a_{l,n}$ into only one set ${\mathcal F}_{j,n}$, to avoid an
overlap of these sets. It is clear that, for any $n \in \N$,
\begin{equation*}
\bigcup_{j=1}^m {\mathcal F}_{j,n} = \{ a_{l,n} \}_{l=1}^n \quad
{\rm and} \quad {\mathcal F}_{k_1,n} \bigcap {\mathcal F}_{k_2,n} =
\emptyset, \ k_1 \neq k_2.
\end{equation*}
We next introduce the factors of $F_n(z)$ by setting
\[
F_{j,n}(z) := \prod_{a_{l,n} \in {\mathcal F}_{j,n}} (z-a_{l,n} ),
\quad j=1,\ldots,m,
\]
so that
\[
\|w^{n_j}F_{j,n}\|_E \ge w^{n_j}(c_j^*) \prod_{a_{l,n} \in {\mathcal F}_{j,n}} |c_j^* - a_{l,n}| = \prod_{a_{l,n} \in {\mathcal F}_{j,n}} d_m^*(a_{l,n}), \quad j=1,\dots,m,
\]
where $n_j:=\deg(F_{j,n}).$ Since the normalized counting measures $\nu_{{\mathcal F}_n}$ in the weighted Fekete points converge to the weighted equilibrium measure $\mu_w$ in the weak* topology, see Theorem 1.3 in \cite[p. 145]{ST}, it follows that
\begin{align*}
\liminf_{n \rightarrow \infty} \left( \prod_{j=1}^m \|w^{n_j} F_{j,n}\|_E \right)^{1/n} &\ge \lim_{n \rightarrow \infty} \left( \prod_{l=1}^n
d_m^*(a_{l,n}) \right)^{1/n} \\ &= \lim_{n \rightarrow \infty} \exp
\left( \frac{1}{n} \sum_{j=1}^n \log d_m^*(a_{k,n} ) \right)  \\
&= \exp \left( \int \log d_m^*(z)\, d\mu_w(z)\, \right),
\end{align*}
because $\log d_m^*(z)$ is continuous in $\C.$ We also have that $\dis\lim_{n \rightarrow \infty} \|w^n F_n\|_E^{1/n} = e^{-F_w}$ by Theorem 1.9 of \cite[p. 150]{ST}, which gives
\[
\liminf_{n \rightarrow \infty} \left( \frac{\prod_{j=1}^m \|w^{n_j} F_{j,n}\|_E}{\|w^n F_n\|_E} \right)^{1/n} \ge e^{C_E(m)}.
\]
To show that $C_E$ cannot be replaced by a smaller constant
independent of $m$, one should essentially repeat the above argument
with $m=n\to\infty.$

\db

\noindent{\it Proof of Theorem \ref{thm3.3}.}

We present a proof for the finite set $E=\{c_k\}_{k=1}^m$ here,
which is sufficient for applications in the proofs of Theorems
\ref{thm3.1} and \ref{3.2}. A proof of the general case will appear
in a separate paper.

Let $M:=\{z\in E: w(z) = \sup_E w = \max_E w\}.$ Our first goal is
to show that $d_E^w(z) = d_M(z)\max_E w$ in a neighborhood of
infinity. Since $E$ is finite, there exists $\varepsilon>0$ such that
\[
w(z) < \max_E w - \varepsilon, \quad z\in E\setminus M.
\]
Suppose that there is a sequence of points $\{z_i\}_{i=1}^{\infty}$
in the plane such that $\lim_{i\to\infty} z_i = \infty$, and the
weighted distance function $d_E^w(z_i)$ is attained at the points of
$E\setminus M$ for each $i\in\N$. It follows that
\[
d_E^w(z_i) = w(t_i)|z_i-t_i| < \left(\max_E w-\varepsilon\right)
\left(|z_i|+\max_{1\le k\le m} |c_k|\right),
\]
where $t_i\in E\setminus M.$ Since $M\subset E$, we have that
$d_M^w(z) \le d_E^w(z),\ z\in\C.$ Hence
\begin{align*}
d_M^w(z_i) &= \max_{t\in M} w(t)|z_i-t| = \max_E w \max_{t\in M}
|z_i-t| = d_M(z_i) \max_E w \\ &\le d_E^w(z_i) < \left(\max_E
w-\varepsilon\right) \left(|z_i|+\max_{1\le k\le m} |c_k|\right).
\end{align*}
If we divide the above inequality by $|z_i|$ and let $|z_i|\to\infty$, then we come to the obvious contradiction $\max_E w \le
\max_E w-\varepsilon.$ Thus there exists $R>0$ such that
\begin{align} \label{5.7}
d_E^w(z) = \max_{t\in M} w(t)|z-t| = d_M(z)\max_E w, \quad |z|>R.
\end{align}

Since $\log w(t) + \log|z-t|$ is a subharmonic function of $z$ in
$\C,$ it follows that
\[
\log d_E^w(z) = \max_{t\in E} \left(\log w(t) +
\log|z-t|\right),\quad z\in\C,
\]
is also subharmonic in the plane, cf. \cite[p. 38]{Ra}. Let
$D_r:=\{z\in\C:|z|<r\},$ and write by the Riesz Decomposition
Theorem \cite[p. 76]{Ra}
\[
\log d_E^w(z) = \int \log|z-t|\,d\sigma_r^w(t) + h_r(z),\quad z\in
D_r,
\]
where $\sigma_r^w$ is a positive Borel measure on $D_r$, and where
$h_r$ is harmonic in $D_r.$ Considering a sequence of disks $D_r$
with $r\to\infty$, we extend $\sigma_r^w$ to the measure
$\sigma_E^w$ on the whole plane. It is known \cite{Boy2,Pr1,LP01}
that
\[
\log d_M(z) = \int \log|z-t|\,d\sigma_M(t),\quad z\in\C,
\]
where $\sigma_M$ is a probability measure on $\C.$ Therefore,
\[
\log d_E^w(z) = \log\max_E w + \log d_M(z) = \log\max_E w + \int
\log|z-t|\,d\sigma_M(t),\quad |z|>R,
\]
by \eqref{5.7}. For any function $u$ that is subharmonic in $\C$,
one can find the Riesz measure of $D_r$ from the formula
\[
\mu(D_r) = r \frac{d}{dr} L(u;r)
\]
except for at most countably many $r$, where
\[
L(u;r):=\frac{1}{2\pi} \int_0^{2\pi} u(re^{i\theta})\,d\theta,
\]
see Theorem 1.2 of \cite[p. 84]{ST}. We remark that Theorem 1.2 is stated in \cite[p. 84]{ST} for potentials of compactly supported measures, but the more general version we use here follows immediately by writing the Riesz decomposition of $u$ on any disk into the sum of a potential and a harmonic function. It is clear that
\[
L(\log d_E^w;r) = L(\log d_M;r) + \log\max_E w,\quad r>R,
\]
so that
\[
\sigma_E^w(D_r) = \sigma_M(D_r),\quad r>R,
\]
except for at most countably many $r$. Consequently, $\sigma_E^w(\C)
= \sigma_M(\C) = 1.$

We also have for any $r>R$ that
\[
\int \log|z-t|\,d\sigma_r^w(t) + h_r(z)= \log\max_E w + \int
\log|z-t|\,d\sigma_M(t),\quad R<|z|<r.
\]
Applying the Unicity Theorem \cite[p. 97]{ST}, we conclude that the two measures coincide in $R<|z|<r$ for any $r>R,$ which gives
\[
\sigma_E^w\vert_{|z|>R} = \sigma_M\vert_{|z|>R}.
\]
This implies that
\[
h_r(z) = \log\max_E w + \int_{|t|\le 2R} \log|z-t|\,d\sigma_M(t) -
\int_{|t|\le 2R} \log|z-t|\,d\sigma_E^w(t),\quad R<|z|<r,
\]
for all $r>R.$ But the right-hand side of this equation is harmonic
and bounded for $|z|>2R$, with the limit value $\log\max_E w$ at
$\infty.$ Thus $h_r$ is continued to a harmonic and bounded function
in $\C$, and it must be identically equal to the constant
$\log\max_E w$ by Liouville's theorem.

\db

\bigskip
\noindent{I. E. Pritsker\\ Department of Mathematics\\ 401
Mathematical Sciences\\ Oklahoma State University\\ Stillwater, OK
74078-1058\\ U.S.A.\\ e-mail: igor@math.okstate.edu}

\bigskip
\noindent{E. B. Saff\\ Center for Constructive Approximation\\
Department of Mathematics\\ Vanderbilt University\\ Nashville,
TN 37240\\ U.S.A.\\
e-mail: Edward.B.Saff@Vanderbilt.Edu}

\end{document}